\documentclass{article}
\usepackage[english,russian,ukrainian]{babel}
\usepackage{latexsym,amssymb,amsmath}
\usepackage[cp1251]{inputenc}
\usepackage{amsfonts,amssymb}

\newcounter{theorem}
\newcounter{lemma}

\renewcommand{\thetheorem}{\arabic{theorem}}

\newcommand{\theor}{\par\refstepcounter{theorem}{\bf Теорема \thetheorem .}\,\,}

\def\Re{\mathop\mathrm{Re}}\,
\def\Im{\mathop\mathrm{Im}}\,

\begin{document}
\Large

\centerline{\textbf{Степеневі ряди та ряди Лорана}}
\centerline{\textbf{в алгебрі комплексних кватерніонів}}
\vskip 5mm
\centerline{\textbf{Т.\,С.~Кузьменко}}
\vskip 5mm

\large{Taylor's and Laurent's expansions of $G$-monogenic mappings taking values in the algebra of complex quaternion are obtained and singularities of these mappings are classified.}

\large{В алгебре комплексных кватернинионов получены тейлоровские и лорановские разложения $G$-моногенных отображений, классифицированы их особые точки.}

\Large
\vskip 5mm

\textbf{1. Вступ.} Нехай $\mathbb{H(C)}$ --- алгебра кватерніонів над полем комплексних чисел
$\mathbb{C}$, базис якої складається з одиниці алгебри $1$ і елементів $I,J,K,$
для яких виконуються наступні правила множення:
$$I^2=J^2=K^2=-1,\,IJ=-JI=K,\,JK=-KJ=I,\,KI=-IK=J.$$

Розглянемо в алгебрі $\mathbb{H(C)}$ інший базис $\{e_1,e_2,e_3,e_4\}$, розклад
 елементів якого в базисі $\{1,I,J,K\}$ має вигляд
$$e_1=\frac{1}{2}(1+iI), \quad e_2=\frac{1}{2}(1-iI), \quad
e_3=\frac{1}{2}(iJ-K), \quad e_4=\frac{1}{2}(iJ+K),$$
де $i$ --- уявна комплексна одиниця. Таблиця множення в новому базисі
набуває вигляду

$$
\begin{tabular}{c||c|c|c|c|}
$\cdot$ & $e_1$ & $e_2$ & $e_3$ & $e_4$\\
\hline
\hline
$e_1$ & $e_1$ & $0$ & $e_3$ & $0$\\
\hline
$e_2$ & $0$ & $e_2$ & $0$ & $e_4$\\
\hline
$e_3$ & $0$ & $e_3$ & $0$ & $e_1$\\
\hline
$e_4$ & $e_4$ & $0$ & $e_2$ & $0$\\
\hline
\end{tabular}\,\,.
$$

Норма кватерніона $c=\sum\limits_{k=1}^4c_ke_k,\,c_k\in\mathbb{C}$
визначається рівністю
\begin{equation}\label{norma-1}
\|c\|:=\sqrt{\sum\limits_{k=1}^4|c_k|^2}\,,
\end{equation}
а
одиниця алгебри $\mathbb{H(C)}$ в цьому базисі є сумою ідемпотентів: $1=e_1+e_2$.

Розглянемо лінійні функціонали $f_1 \colon\mathbb{H(C)\longrightarrow\mathbb{C}}$ та $f_2 \colon\mathbb{H(C)\longrightarrow\mathbb{C}},$ покладаючи
$$f_1(e_1)=f_1(e_3)=1, \qquad f_1(e_2)=f_1(e_4)=0,$$
$$f_2(e_2)=f_2(e_4)=1, \qquad f_2(e_1)=f_2(e_3)=0.$$

Нехай
$$i_1=1, \qquad i_2=a_1e_1+a_2e_2, \qquad i_3=b_1e_1+b_2e_2$$
при $a_k,b_k\in\mathbb{C},\,\,
k=1,2$ --- трійка лінійно незалежних векторів над полем дійсних чисел $\mathbb{R}$ (див. \cite[с.~223]{Plaksa_Puht}).
Це означає, що рівність $$\alpha_1i_1+\alpha_2i_2+\alpha_3i_3=0, \qquad
\alpha_1,\alpha_2,\alpha_3\in\mathbb{R}$$ виконується тоді і тільки тоді,
коли  $\alpha_1=\alpha_2=\alpha_3=0$.

Виділимо в алгебрі $\mathbb{H(C)}$ лінійну оболонку $E_3 \colon =\{\zeta=xi_1+yi_2+zi_3 \colon x,y,z\in\mathbb{R}\}$ над полем дійсних чисел
$\mathbb{R},$ породжену векторами $i_1,i_2,i_3.$ Області $\Omega$
тривимірного простору $\mathbb{R}^3$ поставимо у відповідність область $\Omega_\zeta \colon=\{\zeta=xi_1+yi_2+zi_3:(x,y,z)\in\Omega\}$ в $E_3$.

Введемо позначення $\xi_1 \colon=f_1(\zeta)=x+ya_1+zb_1,\, \xi_2 \colon=
f_2(\zeta)=x+ya_2+zb_2.$ Тепер елемент $\zeta\in E_3$ можна подати у
 вигляді $\zeta=\xi_1e_1+\xi_2e_2$, і згідно визначення (\ref{norma-1})
\begin{equation}\label{norma}
\|\zeta\|=\sqrt{|\xi_1|^2+|\xi_2|^2}.
\end{equation}

Позначимо через $f_k(E_3) \colon=\{f_k(\zeta):\zeta\in E_3\},\, k=1,2$. Відмітимо, що в подальшому  істотним є припущення:
 $f_1(E_3)=f_2(E_3)=\mathbb{C}$.
Очевидно, що воно має місце тоді і тільки тоді, коли
хоча б одне з чисел у кожній з пар $(a_1,b_1)$, $(a_2,b_2)$ належить
 $\mathbb{C}\setminus\mathbb{R}$.

Неперервне відображення $\Phi:\Omega_\zeta\longrightarrow\mathbb{H(C)}$ \big(або $\widehat{\Phi}:\Omega_\zeta\longrightarrow\mathbb{H(C)}$\big) називається
\emph{право-$G$-моногенним}
\big(або \emph{ліво-$G$-моногенним}\big) в області
$\Omega_\zeta\subset E_3$, якщо $\Phi$ \big(або $\widehat{\Phi}$\big)
диференційовне за Гато у кожній точці цієї області, тобто якщо для
кожного $\zeta\in\Omega_\zeta$ існує елемент
$\Phi'(\zeta)$ \big(або $\widehat{\Phi}'(\zeta)$\big) алгебри
$\mathbb{H(C)}$ такий,
що виконується рівність
$$\lim\limits_{\varepsilon\rightarrow 0+0}\Big(\Phi(\zeta+\varepsilon h)-\Phi(\zeta)\Big)\varepsilon^{-1}= h\Phi'(\zeta)\quad\forall\,h\in E_3$$
$$\Biggr(\mbox{або}\,\, \lim\limits_{\varepsilon\rightarrow 0+0}
\left(\widehat{\Phi}(\zeta+\varepsilon h)-\widehat{\Phi}(\zeta)\right)
\varepsilon^{-1}= \widehat{\Phi}'(\zeta)h\quad\forall\,h\in E_3\Biggr).$$

Із розкладу резольвенти \big(див. \cite{Shpakivskiy-Kuzmenko}\big)
$$(t-\zeta)^{-1}=\frac{1}{t-\xi_1}\,e_1+\frac{1}{t-\xi_2}\,e_2\,,\qquad
\forall\,\,t\in \mathbb{C}:\,\,t\neq \xi_1\,, \,\, t\neq \xi_2$$
випливає, що точки $(x,y,z)\in\mathbb{R}^3,$ які відповідають необоротним
елементам $\zeta=xi_1+yi_2+zi_3\in E_3$, лежать на прямих
$$L^1 \colon x+y\Re\,a_1+z\Re\,b_1=0, \qquad y\Im\,a_1+z\Im\,b_1=0,$$
$$L^2 \colon x+y\Re\,a_2+z\Re\,b_2=0, \qquad y\Im\,a_2+z\Im\,b_2=0$$
в просторі $\mathbb{R}^3$. Прямим $L^1,\,L^2$ поставимо у відповідність конгруентні їм прямі $L^1_\zeta\,,\;L^2_\zeta$ в $E_3$.

В роботах \cite{Moisil} --- \cite{Gentili} для різних класів функцій кватерніонної змінної встановлено аналоги класичних результатів комплексного аналізу, зокрема, доведено аналоги інтегральних теорем Коші, інтегральної формули Коші, теореми Морера та досліджено розклади в степеневі ряди і ряди Лорана.

В роботі \cite{Shpakiv-Kuzm-intteor} для $G$-моногенних відображень отримано аналоги інтегральної теореми Коші для криволінійного і поверхневого інтегралів та аналог інтегральної формули Коші.

В цій роботі отримано розклади $G$-моногенних відображень в ряди Тейлора і Лорана, а також здійснено класифікацію особливих точок згаданих відображень.

\textbf{2. Степеневі ряди і $G$-моногенні відображення.}
Безпосереднє застосування способу розкладу голоморфних функцій, що базується на розкладі в ряд ядра Коші (див., наприклад,
\cite[с. 107]{Shabat}), до право-$G$-моногенного в області $\Omega_{\zeta}$ відображення
$\Phi \colon\Omega_{\zeta}\longrightarrow\mathbb{H(C)}$,
дозволяє отримати його розклад в степеневий ряд
\begin{equation}\label{rjad-1}
\Phi(\zeta)=\sum_{n=0}^{\infty} (\zeta-\zeta_0)^{n} p_n
\end{equation}
в кулі з центром у фіксованій точці $\zeta_0\in\Omega_{\zeta}$
меншого радіуса, ніж відстань від точки $\zeta_0$ до межі
$\partial\Omega_{\zeta}$; тут
$$p_n=\frac{\Phi^{(n)}(\zeta_0)}{n!}=\frac{1}{2\pi
i}\int\limits_{\gamma_{\zeta}}\Bigl((\tau-\zeta_0)^{-1}\Bigr)^{n+1}\,d\tau\Phi(\tau)\,,\quad
n=0,1,\dots,$$

а $\gamma_{\zeta}$ визначена в роботі \cite{Shpakiv-Kuzm-intteor}. Це пов'язано з тим, що в
нерівності $\|ab\|\leq 2\sqrt{2}\,\|a\|\|b\|$ константа
$2\sqrt{2}$ не може бути замінена одиницею.

Проте, використовуючи представлення право-$G$-моногенного відображення
$\Phi \colon\Omega_{\zeta}\longrightarrow\mathbb{H(C)}$ через аналітичні
функції комплексної змінної, покажемо, що право-$G$-моногенне відображення $\Phi$ в кулі
$B(\zeta_0,R) \colon=\{\zeta\in E_3:
\|\zeta-\zeta_0\|<R\}$ радіуса
$R \colon=\min\limits_{\zeta\in\partial\Omega_\zeta}\|\zeta-\zeta_0\|$
з центром в довільній точці $\zeta_{0}\in\Omega_{\zeta}$
подається у вигляді ряду (\ref{rjad-1}).

Розглядаючи питання про розклад $G$-моногенних відображень
$\Phi \colon\Omega_\zeta\longrightarrow\mathbb{H(C)}$ і $\widehat{\Phi} \colon\Omega_\zeta\longrightarrow\mathbb{H(C)}$ в степеневий ряд, без
обмеження загальності будемо вважати, що область $\Omega_\zeta\subset
E_3$ є обмеженою.

Нехай $\zeta_{0} \colon=x_{0}i_1+y_{0}i_{2}+z_{0}i_{3}$ --- довільна фіксована точка області $\Omega_\zeta$ і
$R_0 \colon=\min\limits_{\zeta\in\partial\Omega_\zeta}\|\zeta-\zeta_0\|$, де
через $\partial\Omega_\zeta$ позначено межу області $\Omega_\zeta$ в $E_3$.
Візьмемо кулю $\Theta(\zeta_0,R_0) \colon=\{\zeta\in E_3 \colon
\|\zeta-\zeta_0\|<R_0\}$ в $E_3$ радіуса $R_0$ з центром в точці $\zeta_0$.

Позначимо через
$\xi_{10} \colon=x_{0}+a_1y_{0}+b_1z_0$ і
$\xi_{20} \colon=x_{0}+a_2y_{0}+b_2z_0$
точки комплексної площини, які відповідають точці $\zeta_{0}$
за формулами $\xi_{10}=f_1(\zeta_{0})$ і $\xi_{20}=f_2(\zeta_{0})$,
а через $\widetilde{D}_1$ і $\widetilde{D}_2$ --- області в
$\mathbb{C}$, на які куля $\Theta(\zeta_0,R_0)$ відображається
відповідно функціоналами $f_1$ і $f_2$.

Нехай $R \colon=\min\Bigl\{ R_0\,,
\min\limits_{\tau_1\in\partial\widetilde{D}_1}|\tau_1-\xi_{10}|\,,
\min\limits_{\tau_2\in\partial\widetilde{D}_2}|\tau_2-\xi_{20}|\Bigr\}$, де через
$\partial\widetilde{D}_1$ і $\partial\widetilde{D}_2$ позначено
відповідно межі областей
$\widetilde{D}_1$ і $\widetilde{D}_2$. Через $U(\xi_{10},R):=\{\xi_1\in
\mathbb{C}  \colon |\xi_1-\xi_{10}|<R\}$ і $U(\xi_{20},R):=\{\xi_2\in
\mathbb{C}  \colon |\xi_2-\xi_{20}|<R\}$ позначимо круги радіуса $R$ в
комплексній площині з центрами відповідно в точках $\xi_{10}$ і
$\xi_{20}$.

Покажемо, що представлення \big(див. \cite{Shpakivskiy-Kuzmenko}\big)
\begin{equation}\label{Phi-r-rozklad}
\Phi(\zeta)=F_1(\xi_1)e_1+F_2(\xi_2)e_2+F_3(\xi_1)e_3+F_4(\xi_2)e_4,
\end{equation}
де $F_{1}$ і $F_3$ --- деякі голоморфні в
$U(\xi_{10},R)$ функції,  а $F_{2}$ і $F_{4}$ --- деякі
голоморфні в $U(\xi_{20},R)$ функції, дозволяє отримати розклад право-$G$-моногенного
відображення в степеневий ряд
\begin{equation}\label{r-Teylor-series}
\Phi(\zeta)=\sum_{n=0}^{\infty}(\zeta-\zeta_0)^{n}p_n
\end{equation}
в області $B(\zeta_0,R):=\{\zeta\in E_3  \colon f_1(\zeta)\in
U(\xi_{10},R) \,, f_2(\zeta)\in U(\xi_{20},R)\}$,
а представлення
\begin{equation}\label{Phi-l-rozklad}
\widehat{\Phi}(\zeta)=F_1(\xi_1)e_1+F_2(\xi_2)e_2+F_3(\xi_2)e_3+F_4(\xi_1)e_4\,,
\end{equation}
де $F_{1}$ і $F_4$ --- деякі голоморфні в
$U(\xi_{10},R)$ функції,  а $F_2$ і $F_3$ --- деякі
голоморфні в $U(\xi_{20},R)$ функції, дозволяє отримати розклад ліво-$G$-моногенного
відображення в степеневий ряд
\begin{equation}\label{l-Teylor-series}
\widehat{\Phi}(\zeta)=\sum_{n=0}^{\infty}\widehat{p}_n(\zeta-\zeta_0)^{n}
\end{equation}
в області $B(\zeta_0,R)$.

Оскільки за побудовою область $B(\zeta_0,R)$ є
опуклою в напрямку прямих $L^1_\zeta$ і $L^2_\zeta$\,, то в ній право-$G$-моногенне відображення $\Phi(\zeta)$ подається у вигляді
\eqref{Phi-r-rozklad}, а ліво-$G$-моногенне відображення $\widehat{\Phi}(\zeta)$ подається у вигляді
\eqref{Phi-l-rozklad}.

\vskip 1mm

\theor\label{teor-Teylor} \textit{Нехай $f_1(E_3)=f_2(E_3)=\mathbb{C}$, відображення $\Phi \colon\Omega_\zeta\longrightarrow\mathbb{H(C)}$ --- право-$G$-моногенне, а $\widehat{\Phi} \colon\Omega_\zeta\longrightarrow\mathbb{H(C)}$ --- ліво-$G$-моногенне в
довільній обмеженій області $\Omega_\zeta\subset E_3$ і
$\zeta_0\in\Omega_\zeta$. Тоді в області $B(\zeta_0,R)$ відображення $\Phi$
подається у вигляді суми абсолютно збіжного степеневого ряду
\eqref{r-Teylor-series}, в якому
\begin{equation}\label{teor-Teylor-p_n}
p_n=a_ne_1+b_ne_2+c_ne_3+d_ne_4\,,
\end{equation}
і $a_n\,, b_n\,, c_n\,, d_n$ --- коефіцієнти рядів Тейлора функцій}
\begin{equation}\label{teor-Teylor-r-F1,F2,F3,F4}
\begin{array}{c}
\displaystyle
F_1(\xi_1)=\sum_{n=0}^{\infty}a_{n}(\xi_1-\xi_{10})^{n},\qquad
F_2(\xi_2)=\sum_{n=0}^{\infty}b_{n}(\xi_2-\xi_{20})^{n},\\[6mm]
\displaystyle
F_3(\xi_1)=\sum_{n=0}^{\infty}c_{n}(\xi_1-\xi_{10})^{n},\qquad
F_4(\xi_2)=\sum_{n=0}^{\infty}d_{n}(\xi_2-\xi_{20})^{n},
\end{array}
\end{equation}
\emph{які містяться у представленні
\eqref{Phi-r-rozklad} відображення $\Phi(\zeta)$ при $\zeta\in
B(\zeta_0,R)$,
а $\widehat{\Phi}$ подається у вигляді суми абсолютно збіжного степеневого ряду
\eqref{l-Teylor-series}, в якому
\begin{equation}\label{teor-Teylor-p_n-}
\widehat{p}_n=\widehat{a}_ne_1+\widehat{b}_ne_2+\widehat{c}_ne_3+\widehat{d}_ne_4
\end{equation}
і $\widehat{a}_n\,, \widehat{b}_n\,, \widehat{c}_n\,, \widehat{d}_n$ --- коефіцієнти рядів Тейлора функцій}
\begin{equation}\label{teor-Teylor-l-F1,F2,F3,F4}
\begin{array}{c}
\displaystyle
\widehat{F}_1(\xi_1)=\sum_{n=0}^{\infty}\widehat{a}_{n}(\xi_1-\xi_{10})^{n},\qquad
\widehat{F}_2(\xi_2)=\sum_{n=0}^{\infty}\widehat{b}_{n}(\xi_2-\xi_{20})^{n},\\[6mm]
\displaystyle
\widehat{F}_3(\xi_2)=\sum_{n=0}^{\infty}\widehat{c}_{n}(\xi_2-\xi_{20})^{n},\qquad
\widehat{F}_4(\xi_1)=\sum_{n=0}^{\infty}\widehat{d}_{n}(\xi_1-\xi_{10})^{n},
\end{array}
\end{equation}
\emph{які містяться у представленні
\eqref{Phi-l-rozklad} відображення $\widehat{\Phi}(\zeta)$ при $\zeta\in
B(\zeta_0,R)$.}
\vskip 1mm

\textbf{\textit{Доведення.}} Оскільки у представленні \eqref{Phi-r-rozklad} функції $F_1,F_3$ голоморфні в крузі
$U(\xi_{10},R)$, а функції $F_2$, $F_4$ голоморфні в крузі
$U(\xi_{20},R)$, то у відповідних кругах ряди
\eqref{teor-Teylor-r-F1,F2,F3,F4} абсолютно збіжні. Тоді при всіх
$\zeta\in B(\zeta_0,R)$ рівність \eqref{Phi-r-rozklad}
набуває вигляду
$$\Phi(\zeta)=\sum\limits_{n=0}^{\infty}a_{n}(\xi_{1}-\xi_{10})^{n}e_1+\sum\limits_{n=0}^{\infty}b_{n}(\xi_{2}-\xi_{20})^{n}e_2+$$
$$+\sum\limits_{n=0}^{\infty}c_{n}(\xi_{1}-\xi_{10})^{n}e_3+\sum\limits_{n=0}^{\infty}d_{n}(\xi_{2}-\xi_{20})^{n}e_4.$$

Звідси, враховуючи співвідношення
\begin{equation}\label{zeta-xi}
\begin{array}{ll}
(\zeta-\zeta_{0})^{n}e_1=(\xi_1-\xi_{10})^{n}e_1,\quad (\zeta-\zeta_{0})^{n}e_2=(\xi_2-\xi_{20})^{n}e_2,&\\[3mm]
(\zeta-\zeta_{0})^{n}e_3=(\xi_1-\xi_{10})^{n}e_3,\quad (\zeta-\zeta_{0})^{n}e_4=(\xi_2-\xi_{20})^{n}e_4&\\
\end{array}
\end{equation}
для всіх $\zeta\in E_3$ і $n=0,1,\dots$, приходимо до розкладу
\eqref{r-Teylor-series}, коефіцієнти якого визначаються рівністю
\eqref{teor-Teylor-p_n}, при цьому ряд \eqref{r-Teylor-series}
абсолютно збігається в області $B(\zeta_0,R)$. Повністю аналогічно доводиться рівність \eqref{l-Teylor-series}. Теорему доведено.

\vskip 1mm

Тепер так, як і для голоморфних функцій комплексної змінної (див.,
наприклад, \cite[с. 118]{Shabat}), встановлюється наступна теорема
єдиності для $G$-моногенних відображень, які визначені в області
$\Omega_\zeta\subset E_3$ і приймають значення в алгебрі $\mathbb{H(C)}$.

\vskip 1mm

\theor \label{teor-edinosti} {\em Якщо два право-$G$-моногенні відображення
$\Phi_1 \colon\Omega_\zeta\longrightarrow\mathbb{H(C)}$,
$\Phi_2 \colon\Omega_\zeta\longrightarrow\mathbb{H(C)}$
в довільній
області $\Omega_\zeta\subset E_3$ співпадають в деякому околі
довільної внутрішньої точки області $\Omega_\zeta$, то вони тотожно
рівні у всій області $\Omega_\zeta$.}

Аналогічна теорема справедлива і для
 ліво-$G$-моногенних відображень.

 Зазначимо, що співпадання відображень $\Phi_1 \colon\Omega_\zeta\longrightarrow\mathbb{H(C)}$ і $\Phi_2 \colon\Omega_\zeta\longrightarrow\mathbb{H(C)}$ на множині точок, яка містить хоча б одну граничну точку області $\Omega_\zeta$\,, є недостатнім для тотожної рівності цих відображень у всій області $\Omega_\zeta$. Так, наприклад, значення $G$-моногенних в $E_3$ відображень $\Phi_1(\zeta)=\zeta  e_3$ і $\Phi_2(\zeta)=0$ співпадають для всіх $\zeta\in L_\zeta^1$, проте не співпадають тотожно.

\vskip 1mm

\textbf{3. Ряди Лорана і класифікація особливих точок $G$-моногенних відображень.} Розглянемо питання про розклад право-$G$-моногенного відображення
$\Phi \colon\mathcal{K}_\zeta\longrightarrow\mathbb{H(C)}$ і ліво-$G$-моногенного відображення
$\widehat{\Phi} \colon\mathcal{K}_\zeta\longrightarrow\mathbb{H(C)}$ в ряд Лорана відносно точки $\zeta_{0} \colon=x_0i_1+y_0i_2+z_0i_3$, вважаючи, що
відображення задане в необмеженій області
$$\mathcal{K}_\zeta \colon=\{\zeta\in E_3 \colon
0\leq r<{|\xi_1-\xi_{10}|}<R\leq\infty\,,0\leq
r<{|\xi_2-\xi_{20}|}<R\leq\infty\}\,.$$

\vskip 1mm

\theor \label{teor-Lorana} \textit{Нехай $f_1(E_3)=f_2(E_3)=\mathbb{C}$. Тоді кожне право-$G$-моногенне відображення $\Phi \colon\mathcal{K}_\zeta\longrightarrow\mathbb{H(C)}$
в області
$\mathcal{K}_\zeta$ подається в ній у вигляді суми абсолютно збіжного ряду
\begin{equation}\label{r-teor-Loran-series}
\Phi(\zeta)=\sum\limits_{n=-\infty}^{\infty}(\zeta-\zeta_0)^{n}p_n\,,
\end{equation}
де коефіцієнти $p_n$ визначаються формулами \eqref{teor-Teylor-p_n}, в яких $a_n, b_{n}, c_{n}, d_{n}$ --- коефіцієнти рядів Лорана функцій
\begin{equation}\label{r-teor-Lorana-F1,F2,F3,F4}
\begin{array}{c}
\displaystyle
F_1(\xi_1)=\sum_{n=-\infty}^{\infty}a_{n}(\xi_1-\xi_{10})^{n},\qquad
F_{2}(\xi_2)=\sum_{n=-\infty}^{\infty}b_{n}(\xi_2-\xi_{20})^{n},\\[6mm]
\displaystyle
F_3(\xi_1)=\sum_{n=-\infty}^{\infty}c_{n}(\xi_1-\xi_{10})^{n},\qquad
F_4(\xi_2)=\sum_{n=-\infty}^{\infty}d_{n}(\xi_2-\xi_{20})^{n},
\end{array}
\end{equation}
які містяться в розкладі
\eqref{Phi-r-rozklad} відображення $\Phi(\zeta)$,
а ліво-$G$-моногенне відображення $\widehat{\Phi}:\mathcal{K}_\zeta\longrightarrow\mathbb{H(C)}$
подається в області
$\mathcal{K}_\zeta$ у вигляді суми абсолютно збіжного ряду
\begin{equation}\label{l-teor-Loran-series}
\widehat{\Phi}(\zeta)=\sum\limits_{n=-\infty}^{\infty}\widehat{p}_n(\zeta-\zeta_0)^{n}\,,
\end{equation}
де коефіцієнти $\widehat{p}_n$ визначаються формулами \eqref{teor-Teylor-p_n-}, в яких $\widehat{a}_n, \widehat{b}_{n}, \widehat{c}_{n}, \widehat{d}_{n}$ --- коефіцієнти рядів Лорана функцій
\begin{equation}\label{l-teor-Lorana-F1,F2,F3,F4}
\begin{array}{c}
\displaystyle
\widehat{F}_1(\xi_1)=\sum_{n=-\infty}^{\infty}\widehat{a}_{n}(\xi_1-\xi_{10})^{n},\qquad
\widehat{F}_{2}(\xi_2)=\sum_{n=-\infty}^{\infty}\widehat{b}_{n}(\xi_2-\xi_{20})^{n},\\[6mm]
\displaystyle
\widehat{F}_3(\xi_2)=\sum_{n=-\infty}^{\infty}\widehat{c}_{n}(\xi_2-\xi_{20})^{n},\qquad
\widehat{F}_4(\xi_1)=\sum_{n=-\infty}^{\infty}\widehat{d}_{n}(\xi_1-\xi_{10})^{n},
\end{array}
\end{equation}
які містяться в розкладі \eqref{Phi-l-rozklad} відображення $\widehat{\Phi}(\zeta)$ при $\zeta\in
\mathcal{K}_\zeta$. При цьому
 $(\zeta-\zeta_0)^{n}$ при $n=-1,-2,\dots$ в рівностях \em (\ref{r-teor-Loran-series}), (\ref{l-teor-Loran-series}) \em визначається рівністю $(\zeta-\zeta_0)^{n}:=\big((\zeta-\zeta_0)^{-1}\big)^{-n}$.
}
\vskip 1mm

\textbf{\textit{Доведення.}} Оскільки в поданні \eqref{Phi-r-rozklad} функції $F_1,\,F_3$ голоморфні в кільці
$\{\xi_1\in\mathbb{C} \colon r<|\xi_1-\xi_{10}|<R\}$, а функції $F_2$,
$F_4$ голоморфні в кільці
$\{\xi_2\in\mathbb{C} \colon r<|\xi_2-\xi_{20}|<R\},$ то ряди
\eqref{r-teor-Lorana-F1,F2,F3,F4} у відповідних кільцях абсолютно
збіжні. Використовуючи розклади \eqref{r-teor-Lorana-F1,F2,F3,F4}, рівність
\eqref{Phi-r-rozklad} набуває вигляду
$$\Phi(\zeta)=\sum\limits_{n=-\infty}^{\infty}a_{n}(\xi_1-\xi_{10})^{n}e_1+\sum\limits_{n=-\infty}^{\infty}b_{n}(\xi_2-\xi_{20})^{n}e_2+$$
$$+\sum\limits_{n=-\infty}^{\infty}c_{n}(\xi_1-\xi_{10})^{n}e_3+\sum\limits_{n=-\infty}^{\infty}d_{n}(\xi_2-\xi_{20})^{n}e_4.$$

Тепер, враховуючи співвідношення \eqref{zeta-xi}, які
виконуються при всіх цілих значеннях $n$, отримуємо розклад
відображення $\Phi$ в абсолютно збіжний в області $\mathcal{K}_\zeta$ ряд
\eqref{r-teor-Loran-series}, коефіцієнти якого визначаються
рівностями \eqref{teor-Teylor-p_n}. Повністю аналогічно доводиться розклад у ряд Лорана ліво-$G$-моногенного відображення \eqref{Phi-l-rozklad}. Теорему доведено.

\vskip 1mm

Сукупність членів ряду Лорана \eqref{r-teor-Loran-series} \big(або \eqref{l-teor-Loran-series}\big) з
невід'ємними степенями називають його {\em правильною частиною}, а
сукупність членів цього ряду з від'ємними степенями --- {\em
головною частиною} ряду Лорана.

Компактифікуємо алгебру $\mathbb{H(C)}$, додаючи до неї нескінченно
віддалену точку $\infty$, до якої прямує кожна послідовність
$w_n \colon=\tau_{1,n}e_1+\tau_{2,n}e_2+\tau_{3,n}e_3+\tau_{4,n}e_4$, де
$\tau_{1,n},\,\tau_{2,n},\,\tau_{3,n},\,\tau_{4,n}\in\mathbb{C}$, у випадку
коли хоча б одна з послідовностей $\tau_{1,n}$, $\tau_{2,n}$,
$\tau_{3,n}$, $\tau_{4,n}$ збігається до нескінченно віддаленої точки розширеної
комплексної площини.

Припустимо тепер, що право-$G$-моногенне відображення
$\Phi \colon\mathcal{K}^0_\zeta\longrightarrow\mathbb{H(C)}$ і ліво-$G$-моногенне відображення
$\widehat{\Phi} \colon\mathcal{K}^0_\zeta\longrightarrow\mathbb{H(C)}$ задані в області
$$\mathcal{K}^0_\zeta \colon=\{\zeta\in E_3 \colon
0<{|\xi_1-\xi_{10}|}<R\leq\infty\,,
0<{|\xi_2-\xi_{20}|}<R\leq\infty\}\,.$$
Позначимо через $\widetilde{\mathcal{K}}^0_\zeta \colon=\{\zeta\in E_3:
{|\xi_1-\xi_{10}|}<R\,, {|\xi_2-\xi_{20}|}<R\}.$

Справедлива наступна теорема.

\vskip 1mm

\theor\label{teor-Loran-clas} \textit{Нехай $f_1(E_3)=f_2(E_3)=\mathbb{C}$.
Якщо розклад \eqref{r-teor-Loran-series} відображення $\Phi:
\mathcal{K}^0_\zeta\longrightarrow\mathbb{H(C)}$ :}

1) {\em не містить головної частини, то відображення $\Phi$ має
скінченні границі}
\begin{equation}\label{teor-Loran-clas-1}
\lim\limits_{\begin{array}{c}
\zeta\longrightarrow\zeta_0+e^*,\\
\zeta\notin\bigl\{\zeta_{0}+e^* \colon e^*\in L_\zeta^1\cup L_\zeta^2\bigr\}
\end{array}} \Phi(\zeta)
\end{equation}

2) {\em містить лише скінченне число  доданків у головній частині,
то хоча б при одному значенні $k=1$ або $k=2$ відображення
$\Phi$ має нескінченні границі
\begin{equation}\label{teor-Loran-clas-2}
\lim\limits_{\begin{array}{c}
\zeta\longrightarrow\zeta_0+e_k^*,\\
\zeta\notin\bigl\{\zeta_{0}+e_k^* \colon e_k^*\in L_\zeta^k\bigr\}
\end{array}} \Phi(\zeta)
\end{equation}
в усіх точках
$\zeta_0+e_k^*\in\widetilde{\mathcal{K}}^0_\zeta\cap\{\zeta_{0}+e^*_k \colon\,\,e^*_k\in L_\zeta^k\}$;}

3) {\em містить нескінченне число доданків у головній частині, то хоча б при одному значенні $k=1,2$ відображення $\Phi$ або має нескінченну границю,
або не має ні скінченної, ні нескінченної границі в усіх точках
$\zeta_0+e_k^*\in\widetilde{\mathcal{K}}^0_\zeta\cap\{\zeta_{0}+e^*_k \colon\,\,e^*_k\in L_\zeta^k\}$.}

{\em  Аналогічні твердження справедливі для
 ліво-$G$-моногенних відображень.}

\vskip 1mm
\textbf{\textit{Доведення.}}   Відображення $\Phi$ в області $\mathcal{K}_\zeta^0$
подається у вигляді \eqref{Phi-r-rozklad}, де функції
$F_{1},\,F_3$ голоморфні в проколотому околі
$U(\xi_{10},R)\setminus\{\xi_{10}\}$ точки $\xi_{10}$,
а функції  $F_{2}$, $F_{4}$ голоморфні в проколотому околі
$U(\xi_{20},R)\setminus\{\xi_{20}\}$ точки
$\xi_{20}$\,.

Розглянемо випадок, коли розклад \eqref{r-teor-Loran-series} не
містить головної частини, тобто має вигляд \eqref{r-Teylor-series}. При
цьому коефіцієнти рядів Лорана \eqref{r-teor-Lorana-F1,F2,F3,F4}
пов'язані з коефіцієнтами ряду \eqref{r-Teylor-series} співвідношеннями
\eqref{teor-Teylor-p_n}, з яких, в силу рівностей $p_n=0$ при
$n=-1,-2,\ldots\,$, слідують рівності $a_n=b_n=c_n=d_n=0$ при всіх
від'ємних індексах $n$. Отже, ряди Лорана
\eqref{r-teor-Lorana-F1,F2,F3,F4} в околах відповідних точок
$\xi_{10}$, $\xi_{20}$ є рядами Тейлора своїх сум і тому функції
$F_1$, $F_2$, $F_3$, $F_4$ з рівності \eqref{Phi-r-rozklad} є
голоморфними у відповідних областях $U(\xi_{10},R)$ чи
$U(\xi_{20},R)$. Тому відображення \eqref{Phi-r-rozklad} має
скінченні границі \eqref{teor-Loran-clas-1} в усіх точках
$\zeta_0+e^*
\in\widetilde{\mathcal{K}}^0_\zeta\cap\bigl\{\zeta_{0}+e^* \colon
e^*\in L_\zeta^1\cup L_\zeta^2\bigr\}$.

Розглянемо тепер випадок, коли головна частина розкладу
\eqref{r-teor-Loran-series} містить лише скінченне число
доданків, тобто в
\eqref{r-teor-Loran-series} тільки скінченне число відмінних
від нуля коефіцієнтів $p_n$ при від'ємних $n$. Тоді із
співвідношень \eqref{teor-Teylor-p_n}, які пов'язують
коефіцієнти рядів Лорана \eqref{r-teor-Lorana-F1,F2,F3,F4} з
коефіцієнтами ряду \eqref{r-teor-Loran-series}, випливає, що
всі головні частини рядів \eqref{r-teor-Lorana-F1,F2,F3,F4} не
містять нескінченної кількості доданків і головна частина хоча б
одного з них відмінна від тотожного нуля. Тому точка $\xi_{10}$ не
є істотно особливою точкою для функцій $F_1,\,F_3$, а точка $\xi_{20}$
---  для функцій $F_2$ і $F_4$, але хоча б одна з функцій $F_1$,
$F_2$, $F_3$, $F_4$ має полюс у відповідній точці.  Звідси випливає, що
хоча б одна з функцій  $F_1$,
$F_2$, $F_3$, $F_4$ має нескінченну границю
при $\xi_1\rightarrow\xi_{10}$ чи при $\xi_2\rightarrow\xi_{20}$,
тобто границя \eqref{teor-Loran-clas-2} є також нескінченною
при $k=1$ або $k=2$.

Розглянемо, нарешті, випадок, коли головна частина розкладу
\eqref{r-teor-Loran-series} містить нескінченно багато
відмінних від нуля членів, тобто існує нескінченно багато
відмінних від нуля коефіцієнтів $p_{n}$ при відємних $n$. Тоді із
співвідношень \eqref{teor-Teylor-p_n} випливає, що головна
частина хоча б одного з рядів \eqref{r-teor-Lorana-F1,F2,F3,F4}
містить нескінченно багато доданків, а це, в свою чергу, означає,
що або точка $\xi_{10}$ є істотно особливою для функцій $F_1$ чи $F_3$, або
точка $\xi_{20}$ є істотно особливою, принаймні, для однієї з
функцій: $F_2$ чи $F_4$. Тому відображення  $\Phi$ не може мати
скінченної границі в усіх точках множини
$\widetilde{\mathcal{K}}_\zeta^0\cap\{\zeta_{0}+e^* \colon
e^*\in L_\zeta^1\cup L_\zeta^2\}$, але вона може мати в цих
точках нескінченну границю.

Так, наприклад, якщо $\xi_{10}$ --- полюс функції  $F_1$ і істотно особлива точка функції $F_3$, а
$\xi_{20}$ --- істотно особлива точка функцій $F_2$, $F_4$, то функція
$F_1$ має нескінченну границю в точці $\xi_{10}$, а отже, границя
\eqref{teor-Loran-clas-2} є нескінченною в усіх
точках $\zeta_0+e^*\in\widetilde{\mathcal{K}}_\zeta^0\cap\{\zeta_{0}+e^* \colon e^*\in L_\zeta^1\}$.

У випадку, коли, наприклад, $F_2\equiv 0$, $F_3\equiv 0$, $F_4\equiv 0$ і точка $\xi_{10}$ є
істотно особливою для функції $F_1$, відображення $\Phi$ не має ні
скінченної, ні нескінченної границі
\eqref{teor-Loran-clas-2} в усіх точках
$\zeta_0+e^*\in\widetilde{\mathcal{K}}_\zeta^0\cap\{\zeta_{0}+e^* \colon e^*\in L_\zeta^1\}$. Теорему доведено.

\vskip 2mm

Поняття усувної особливої точки, полюса або істотно особливої
точки для $G$-моногенного в області $\mathcal{K}_\zeta^0$ відображення $\Phi \colon
\mathcal{K}_\zeta^0 \longrightarrow\mathbb{H(C)}$ вводяться так, як і
відповідні поняття для голоморфних функцій в комплексній площині
(див., наприклад, \cite[с. 135]{Shabat}). А саме, точка $\zeta_0$
називається:

 1)  \textit{усувною особливою точкою} відображення $\Phi$, якщо існує скінченна границя
$$\lim\limits_{\begin{array}{c}
\zeta\longrightarrow\zeta_{0},\\
\zeta\notin\{\zeta_{0}+e^* \colon e^*\in L_\zeta^1\cup L_\zeta^2\}
\end{array}} \Phi(\zeta)=A;
$$

2)  \textit{полюсом} функції $\Phi$,  якщо існує нескінченна
границя
$$\lim\limits_{\begin{array}{c}
\zeta\longrightarrow\zeta_{0},\\
\zeta\notin\{\zeta_{0}+e^* \colon e^*\in L_\zeta^1\cup L_\zeta^2\}
\end{array}}\Phi(\zeta)=\infty;$$

3)  \textit{істотно особливою точкою} відображення $\Phi$, якщо відображення
$\Phi$ не має ні скінченної, ні нескінченної границі при
$\zeta\longrightarrow\zeta_{0}$ і $\zeta\notin\{\zeta_{0}+e^* \colon e^*\in L_\zeta^1\cup L_\zeta^2\}$.

Повністю аналогічно дані поняття вводяться і для ліво-$G$-моногенних відображень $\widehat{\Phi} \colon
\mathcal{K}^0_\zeta \longrightarrow\mathbb{H(C)}$.

З теореми \ref{teor-Loran-clas} випливає, що ізольована
особлива точка у $G$-моногенного відображення може бути лише усувною, а у
випадку, коли відображення має неусувну особливість в точці
$\zeta_0$, особливими є всі точки хоча б однієї з множин
$\widetilde{\mathcal{K}}^0_\zeta\cap\{\zeta_{0}+e_1^* \colon\,\,e_1^*\in
L_\zeta^1\}$ або $\widetilde{\mathcal{K}}^0_\zeta\cap\{\zeta_{0}+e_2^* \colon\,\,e_2^*\in
L_\zeta^2\}$.

\end{document}